\documentclass[a4paper]{amsart}

\usepackage{amsmath,amsfonts,amssymb,bm}

\usepackage{cancel}
\usepackage{tikz-cd}
\usepackage{float}
\usepackage[all]{xy}
\usepackage{ulem}
\def\lcf{\lbrack\! \lbrack}
\def\rcf{\rbrack\! \rbrack}
\usepackage[latin1]{inputenc}

\usepackage[numbers,sort&compress]{natbib}
\usepackage{graphicx}

\usepackage{hyperref}
\hypersetup{
   colorlinks   =  true,
    linkcolor    = cyan,
    citecolor    = red,
     urlcolor	=magenta,     
}

\usepackage{amsthm}

\newtheorem{definition}{Definition}[section]

\newtheorem{theorem}[definition]{Theorem}
\newtheorem{proposition}[definition]{Proposition}
\newtheorem{corollary}[definition]{Corollary}
\newtheorem{remark}[definition]{Remark}
\newtheorem{example}[definition]{Example}

\newenvironment{pf}{\begin{trivlist}\item[]{\sc Proof.}}%
           {\nolinebreak \hfill $\Box$ \end{trivlist}}

\def \dif {\mathrm{d}}
\def \Hs{\mathrm{Hess}\,}
\def \g{\mathfrak{g}}
\def \R{\mathbb{R}}
\def \C{\mathbb{C}}
\def \Lie{\mathcal{L}}
\newcommand{\parcial}[1]{\frac{\partial}{\partial #1}}
\def \p {\eta}

\newcommand{\cbl}{\color{blue}}

\DeclareMathOperator{\Tr}{Tr}
\DeclareMathOperator{\sgn}{sgn}

\title[Unimodularity and Hamiltonian dynamics on Poisson-Lie groups]{Unimodularity and invariant volume forms for Hamiltonian dynamics on Poisson-Lie groups}
\author{I. Gutierrez-Sagredo$^1$, D. Iglesias Ponte$^2$, J. C. Marrero$^2$, E. Padr\'on$^2$, Z. Ravanpak$^3$}
\thanks{AMS Mathematics Subject Classification (2020): 37C40, 37J39, 53D17, 70G45, 70H05}
\thanks{Keywords: Modular vector fields, modular class, invariant volume forms, Hamiltonian systems, Poisson-Lie groups}

\begin{document}

\maketitle

\vspace{-20pt}
\begin{center}
{\small\it  $\;^1$ Departamento de Matem\'aticas y Computaci\'on, Universidad de Burgos, 09001 Burgos, Spain}

{\small\it e-mail:  igsagredo@ubu.es}

{\small\it $\;^2$ULL-CSIC Geometr\'{\i}a Diferencial y Mec\'anica Geom\'etrica, Departamento de Matem\'aticas, Estad\'{\i}stica e Investigaci\'on Operativa
and Instituto de Matem\'aticas y Aplicaciones (IMAULL)}\\{\small\it University of La Laguna, Spain}
\\[5pt]
{\small\it e-mail: diglesia@ull.edu.es, jcmarrer@ull.edu.es, mepadron@ull.edu.es}

{\small\it $\;^3$ Institute of Mathematics, 
Polish Academy of Sciences, Warsaw, Poland}

{\small\it e-mail: zravanpak@impan.pl}
\end{center}

\begin{abstract} In this paper, we discuss several relations between the existence of invariant volume forms for Hamiltonian systems on Poisson-Lie groups and the unimodularity of the Poisson-Lie structure. In particular, we prove that Hamiltonian vector fields on a Lie group endowed with a unimodular Poisson-Lie structure preserve a multiple of any left-invariant volume on the group. Conversely, we also prove that if there exists a Hamiltonian function such that the identity element of the Lie group is a nondegenerate singularity and the associated Hamiltonian vector field preserves a volume form, then the Poisson-Lie structure is necessarily unimodular. Furthermore, we illustrate our theory with different interesting examples, both on semisimple and unimodular Poisson-Lie groups.

\end{abstract}

\tableofcontents

\section{Introduction}
\label{sec:intro}
\subsection{Integrability for dynamical systems and invariant volume forms}
Given a dynamical system $X\in \mathfrak{X}(M)$ on a manifold $M$ of dimension $n$, the exact integrability is closely related to finding tensors which are invariant for the system. It is clear that if one finds a family of independent $C^\infty$-functions $\{ f_1 ,\ldots ,f_{n-1}\}$ which are first integrals for $X$, i.e., $X(f_i)=0$, then it is possible to obtain, up to reparametrization, the orbits for $X$.

Another approach to the integrability of the dynamical system on $M$, following Euler and Jacobi, is to look for $n-2$ first integrals $\{f_1,\ldots ,f_{n-2}\}$ which are functionally independent and an invariant volume form $\Omega$. Under these conditions, the system can be integrated by quadratures, that is, 
to determine the trajectories of $X$ by means of a finite number of algebraic operations and quadratures of some functions. For instance, if $dim\, M=2$, an invariant volume form $\Omega$ with respect to $X$ implies that
$i_X\Omega$ is a closed 1-form, so locally can be written as $df$, $f$ being  a first integral of $X$ (for more details see, for instance, \cite{Ko} and the references therein).

\subsection{Invariant volume forms for Hamiltonian systems with respect to (almost) Poisson structures}
Among dynamical systems, those that are Hamiltonian with respect to a symplectic structure are of particular interest. In this direction, it is well-known that a Hamiltonian system in a symplectic manifold of dimension $2n$ is completely integrable if it admits $n$ functionally independent first integrals in involution with respect to the Poisson bracket. This kind of systems can be integrated by quadratures. On the other hand, with respect to invariant volume forms, we have Liouville's theorem:

{\it Given a Hamiltonian  on a symplectic manifold, the flow of the Hamiltonian vector field preserves the symplectic volume} (see, for instance, \cite{AM,Ar}).
  
In the more general case when the dynamical system is Hamiltonian with respect to a (not necessarily symplectic) Poisson structure, there are two possible approaches in order to find invariant volume forms: the first one is to describe the symplectic leaves of the Poisson structure and apply Liouville's theorem obtaining invariant volume forms on the leaves; the second is to look for an invariant volume form on the whole manifold. A disadvantage that shows up for the first approach is that the symplectic leaves (and the corresponding symplectic structures on them) for some types of Poisson manifolds are hard to compute. Therefore, in these cases it seems more natural to use the second approach, finding invariant forms on the whole manifold. In this direction, we have the following result:

{\it Given a linear Poisson structure on a vector bundle $A\to M$ and a Hamiltonian function of 
mechanical type on $A$, the Hamiltonian dynamics preserves a volume form of basic type 
if and only if the Poisson structure on $A$ is unimodular. In fact, if the Hamiltonian function is 
of kinetic type then the Hamiltonian dynamics preserves a volume form (not necessarily of basic 
type) if and only if the Poisson structure on $A$ is unimodular } (for more details, see \cite{M}).

We recall that an orientable Poisson manifold $M$ is said to be unimodular if there exists a volume form invariant for all Hamiltonian vector fields. Equivalently, unimodularity can be expressed as the vanishing of the so-called modular class.
The modular class is a first order cohomology class in the Poisson 
cohomology complex of $M$ defined by a Poisson vector field, which measures the 
existence of an invariant volume form for the flows of all Hamiltonian vector fields (see 
\cite{W96}). So, by Liouville's theorem, the modular class of a symplectic manifold is zero. Besides the importance of the 
modular class in a Poisson manifold for dynamics, unimodularity has also been used to prove the 
existence of a duality between Poisson homology and Poisson cohomology \cite{Xu}.

Note that in the linear Poisson case, the existence of an invariant volume form for certain kinds of Hamiltonian systems is closely related to the unimodularity of the Poisson structure. This is still true for the more general class of reduced symmetric nonholonomic  mechanical systems subjected to linear nonholonomic constraints. In fact, there is a geometric formulation for these systems in terms of almost Poisson structures (see, for instance, \cite{LeMaMa,GrLeMaMa}) which is used, in \cite{FGM}, to discuss the existence of invariant volume forms in terms of the unimodularity of such almost Poisson structures.
{
\subsection{Hamiltonian systems on Poisson-Lie groups}
 A Poisson-Lie group is a Lie group $G$ endowed with a Poisson structure such that the multiplication is a Poisson map (see \cite{Dr}). If the Lie group $G$ is a vector space $V$, then $V=\mathfrak{g}^*$ is the dual space of a Lie algebra $\mathfrak{g}$, endowed with the Lie-Poisson structure. Poisson-Lie groups are in one-to-one correspondence with the so-called Lie bialgebras $(\mathfrak{g},\mathfrak{g}^*)$, which are compatible pairs of Lie algebras in duality (for more details, see \cite{CFM,Vaisman} or the recent survey paper \cite{Me} and the references therein). 
 
We remark that some Lie bialgebras can be obtained from solutions of the classical Yang-Baxter equation (on a Lie algebra) and these solutions are specially relevant in the theory of integrable systems. This gives an interesting connection between integrable systems and Poisson-Lie groups, particularly in the case of  factorisable Poisson-Lie groups as it is shown, for instance, in \cite{STS,Se-Tian-Shansky}.

\noindent On the other hand, Hamiltonian systems on Poisson-Lie groups appear in the differential equation approach to the singular value decomposition (SVD) of a bidiagonal matrix \cite{CHU}.  In fact, in \cite{Tomei}, the authors showed that the system considered in \cite{CHU} is Hamiltonian with respect to the (standard) Sklyanin bracket $\{\cdot, \cdot\}$ defined on $\mathrm{SL}(n,\mathbb R)$.

The previous facts provide a good motivation to discuss Hamiltonian systems on Poisson-Lie groups. In particular, the existence of invariant volume forms for such systems. We note that the description of the symplectic leaves for a Poisson-Lie group and the symplectic structure on them may be hard. In addition, one could consider the symplectic groupoid ${\mathcal G}(G)$ integrating the Poisson-Lie group $G$ (see \cite{Lu}) and the invariant induced Hamiltonian system on ${\mathcal G}(G)$. Then, one could apply Liouville theorem to  the symplectic Hamiltonian system. The problem is that $\dim {\mathcal G}(G)=2\dim G$ and, therefore, we would  need to duplicate the number of variables. So, for these previous reasons, it seems reasonable to discuss the existence of invariant volume forms directly in the Poisson-Lie group. In this direction, a first result which can be found in the literature has been obtained in \cite{K88} (see also \cite{M}): The flow of a Hamiltonian function of kinetic type on the dual space $\mathfrak{g}^*$ of a Lie algebra $\mathfrak{g}$ (endowed with the Lie-Poisson structure) preserves a volume form if and only if the Lie algebra $\mathfrak{g}$ is unimodular or, equivalently, the Lie-Poisson structure on $\frak{g}^*$ is unimodular. 
}

\subsection{Aim of the paper}
The aim of the paper is to discuss the relation between the existence of an invariant volume form for a Hamiltonian system on a Poisson-Lie group $G$ and the unimodularity of the Lie-Poisson structure on $G$. In fact, the main result of the paper is (see Corollary \ref{first-theorem} and Theorem \ref{th3-2}):

{\it Given a Hamiltonian function $H$ on a Lie group $G$ endowed with a unimodular Poisson-Lie structure $\Pi$ (that is, the dual Lie algebra of $G$ is unimodular), then the Hamiltonian vector field of $H$, $X_{H}^\Pi$, preserves a multiple of any left-invariant volume form on $G$. Conversely, if the identity element of $G$ is a nondegenerate singularity of $H$ and $X_{H}^\Pi$ preserves a volume form then the Poisson-Lie structure $\Pi$ is unimodular.}    

We remark that the modular vector field of a Poisson-Lie group $G$ with respect to a left-invariant volume form has been described in \cite{ELW}, using the modular character of the dual Lie algebra of $G$.

Apart from the previous results mentioned above, we also obtain an interesting necessary condition for a Hamiltonian system on a Poisson-Lie group to preserve a volume form. This condition must be satisfied in the singular points of the Hamiltonian dynamics (see Remark \ref{r3.1'}).

\subsection{Organization of the paper} The paper is organized as follows. In Section 2, we recall the definition of the Poisson cohomology associated with any Poisson structure, the basic properties of Poisson-Lie groups and the definition of the modular vector field associated with any Poisson structure.
In Section 3, we prove the main results of the paper: Corollary \ref{first-theorem} and Theorem \ref{th3-2} (see also Remark \ref{r3.1'}). In Section 4, we describe several examples which illustrate our theory:
\begin{itemize}
\item
We deduce, as a consequence, the previous result by Kozlov \cite{K88} on the existence of invariant volume forms for Hamiltonian systems of kinetic type with respect to Lie-Poisson structures on the dual space of a Lie algebra.
\item
We present an interesting connection between information geometry and Hamiltonian functions on a Poisson-Lie group with a nondegenerate singularity at the identity element.
\item
We discuss the existence of invariant volume forms for two Hamiltonian systems on two semisimple Poisson-Lie groups:
$\mathrm{SL}(2,\mathbb{R})$ and $S^3 \simeq \mathrm{SU}(2,\mathbb{C})$. The Hamiltonian system on $\mathrm{SL}(2,\mathbb{R})$ is just the  dynamical system associated with the SDV problem for a bidiagonal matrix (see \cite{CHU,Tomei}). 
\item
Finally, we present some interesting examples of Hamiltonian systems on Poisson-Lie groups endowed with unimodular Poisson-Lie structures (so, all of them admit invariant volume forms).

\end{itemize}

\section{Poisson structures, Poisson-Lie groups and modular vector fields}
In this section, we will review some definitions and basic constructions on Poisson geometry (for more details, see \cite{CFM,Me,Vaisman}). 

\subsection{Poisson cohomology and symplectic foliation}
A $2$-vector $\Pi$ on a smooth manifold $M$ is said to be {\it Poisson } if 
\begin{equation}\label{Poisson-condition}
\lcf \Pi,\Pi\rcf=0,
\end{equation}
where $\lcf \cdot,\cdot\rcf$ is the Schouten-Nijenhuis bracket on $M$, which is a graded Lie bracket on the space of multivector fields extending the Lie bracket of vector fields. In such a case, $\Pi$ is a {\it Poisson structure} on $M.$ Equivalently, {a Poisson structure} on $M$ is a bracket of real $C^\infty$-functions on $M$, 
\[
\{\cdot,\cdot\}:C^\infty(M)\times C^\infty(M)\to C^\infty(M),
\]
which is ${\mathbb{R}}$-bilinear, skew-symmetric, a derivation in each argument with respect to the standard product on $C^\infty(M)$ and, in addition, it satisfies the Jacobi identity. The relation between the two definitions of a Poisson structure is 
\[
\{F,H\}=\Pi(\dif F,\dif H),\mbox{ for }F,H\in C^\infty(M).
\]

The Poisson $2$-vector $\Pi$ induces, in a natural way, a morphism of vector bundles (over the identity of $M$)
\[
\Pi^\#:T^*M\to TM
\]
given by 
\[
\langle \Pi^\#(\alpha),\beta\rangle =\Pi(x)(\alpha,\beta),\;\;\mbox{ for $x\in M$ and $\alpha,\beta\in T^*_xM.$}
\]
The Hamiltonian vector field $X_H^\Pi$ associated with a real $C^\infty$-function on  $M$, $H:M\to \R,$  is just 
$$X_H^\Pi(x)=\Pi^\#(\dif H(x)), \mbox{ for } x\in M.$$ 
On the other hand, using (\ref{Poisson-condition}), one may introduce a cohomology operator $\partial_\Pi$ on the space of multivector fields on $M.$ Indeed, if $P$ is a $k$-vector on $M$ then $\partial_\Pi P$ is the $(k+1)$-vector defined by 
$$\partial_\Pi P=\lcf\Pi, P\rcf.$$ 
So, $\partial_\Pi H$ is just the Hamiltonian vector field of $H.$

The resultant cohomology is called the {\it Poisson cohomology} of $M.$ Note that the first cohomology group is the quotient 
$$H_\Pi^1(M)=\frac{P_\Pi (M)}{Ham_{\Pi}(M)},$$
where $Ham_{\Pi}(M)$ is the space of Hamiltonian vector fields on $M$ and $P_\Pi(M)$ is the space of Poisson vector fields, that is, the set of vector fields $X$ on $M$ satisfying 
\begin{equation}\label{eq:Poisson:vector:field}
\lcf X,\Pi\rcf={\mathcal L}_X\Pi=0.
\end{equation}
In the particular case when the morphism $\Pi^\#:T^*M\to TM$ is an isomorphim, we can consider the $2$-form $\Omega$ on $M$ given by 
$$\Omega(x)(u,v)=\langle (\Pi^\#)^{-1}(u),v\rangle $$
for $x\in M$ and $u,v\in T_xM.$ The $2$-form  $\Omega$  is a {\it symplectic structure} on $M$, that is, $\Omega$ is closed and non-degenerate. 

In the general case, we have a generalized foliation ${\mathcal F}$ on $M,$ whose characteristic space at the point $x\in M$ is 
$${\mathcal F}(x)=\Pi^\#(T_x^*M),$$
and the Poisson structure induces a symplectic structure on the leaves of this foliation.  $\mathcal F$ is called {\it the symplectic foliation of $M$.}

\subsection{ Poisson-Lie groups}
Poisson-Lie groups are important examples of Poisson manifolds.

A Poisson structure $\Pi$ on a Lie group $G$ is said to be {\it multiplicative } if the multiplication $m:G\times G\to G$ is a Poisson map, when on the product manifold $G\times G$ we consider the standard product Poisson structure, that is,
\[
\Pi (gh)=(L_g)_\ast (\Pi (h))+(R_h)_\ast (\Pi (g))
\]
for $g,h\in G$. In these conditions, $\Pi$ is a {\it Poisson-Lie structure} on $G$. 

In such a case, we have that $\Pi(e)=0$, with $e$ the identity element in $G$. Moreover, if ${\mathfrak g}$ is the Lie algebra of $G$, the differential of $\Pi$ at the identity $e$, $\dif _e\Pi$, induces an adjoint $1$-cocycle 
\[
\delta=\dif _e\Pi\colon{\mathfrak g}\to \wedge^2 {\mathfrak g}.
\]
This means that 
\[
\delta[\xi,\eta]_{\mathfrak g}=[\xi,\delta\eta]_{\mathfrak g}-[\eta,\delta \xi]_{\mathfrak g},
\]
for all $\xi,\eta\in {\mathfrak g}$, where $[\cdot,\cdot]_{\mathfrak g}$ also denotes the algebraic Schouten bracket, that is, the natural extension of the Lie bracket $[\cdot,\cdot]_{\mathfrak g}$ on ${\mathfrak g}$ making on $\oplus_k \wedge^k \mathfrak{g}$ a graded Lie algebra. In addition, \eqref{Poisson-condition} implies that the dual morphism of $\delta$ 
$$\delta^*=[\cdot,\cdot]_{\mathfrak g^*}:\wedge^2{\mathfrak g}^*\to {\mathfrak g}^*$$
is a Lie bracket on ${\mathfrak g}^*.$ 

In other words, the pair $(({\mathfrak g},[\cdot,\cdot]_\g),({\mathfrak g}^*,[\cdot,\cdot]_ {\mathfrak g^*}))$ is {\it a Lie bialgebra.} Conversely, if $(({\mathfrak g},[\cdot,\cdot]_\g,({\mathfrak g}^*,[\cdot,\cdot]_ {\mathfrak g^*}))$ is a Lie bialgebra and $G$ is a connected and simply-connected Lie group with Lie algebra $({\mathfrak g},[\cdot,\cdot]_\g)$ then there exists a unique Poisson-Lie structure $\Pi$ on $G$ such that the Lie bialgebra associated with
$\Pi$ is just 
\[
(({\mathfrak g},[\cdot,\cdot]_\g),({\mathfrak g}^*,[\cdot,\cdot]_ {\mathfrak g^*})).
\]

A first example of a non-trivial Poisson-Lie structure is the linear Poisson structure $\Pi_{LP}$  on the dual space ${\mathfrak g}^*$ of a Lie algebra $({\mathfrak g},[\cdot,\cdot]_\g),$ which is characterized as follows. If $F,G\in C^\infty({\mathfrak g}^*)$  then 
\begin{equation}\label{Lie-Poisson}
\Pi_{LP}(\mu)(\dif F(\mu),\dif G(\mu))=\mu[\dif F(\mu),\dif G(\mu)]_\g, \mbox{ for }\mu\in {\mathfrak g}^*.
\end{equation}

Note that $\dif F(\mu),\dif G(\mu)\in T_\mu^*{\mathfrak g}^*\cong {\mathfrak g}.$ The couple $({\mathfrak g}^*,\Pi_{LP})$ is an abelian  Poisson-Lie group. 

Other interesting examples of Poisson-Lie groups may be obtained from solutions of the generalized Yang-Baxter equation on a Lie algebra or $r$-matrices.  More precisely, if $G$ is a Lie group with Lie algebra $\g$ and $r\in \wedge^2\g$ is a solution of the generalized Yang-Baxter equation, that is, $[r,r]_\g$ is ad-invariant, then 
$$\Pi=r^l-r^r$$ is a Poisson-Lie structure on $G$. Here, if $P\in \wedge^k\g$ then $P^l$ (resp. $P^r$) is the left-invariant (resp. right-invariant) $k$-vector field on $G$ whose value at the identity element of $G$ is $P$. Note that the adjoint $1$-cocycle $\delta: \g\to \wedge^2\g$ associated with a $r$-matrix is given by $$ \delta(\xi)=\mbox{ ad}_\xi r,\mbox{ for } \xi\in \g.$$ More examples of (non-abelian) Poisson-Lie groups will appear throughout the paper. 

\subsection{Modular vector field of a Poisson manifold}


Let $(M,\Pi )$ be an orientable Poisson manifold. Given a volume form $\Phi$, we consider the vector field $\mathcal{M}_\Phi ^\Pi\in{\mathfrak X}(M)$ given by
\begin{equation}\label{eq:modular:vf}
\mathcal{M}_\Phi ^\Pi (H)=\mbox{div}_\Phi (X^\Pi _H), \mbox{ for }H\in C^\infty (M),
\end{equation}
where $\mbox{div}_\Phi (X^\Pi_H)$ is the divergence of the Hamiltonian vector field $X^\Pi_H$ with respect to the volume form $\Phi,$ that is,
\begin{equation}\label{eq:divergencia}
\Lie _{X^\Pi _H} \Phi = \mbox{div}_\Phi (X^\Pi _H) \Phi.
\end{equation}
The vector field $\mathcal{M}_\Phi ^\Pi$ is the {\it modular vector field} of $(M,\Pi)$ with respect to $\Phi$. It is proved in \cite{W96} that $\mathcal{M}_\Phi ^\Pi$ is a Poisson vector field (see \eqref{eq:Poisson:vector:field}).
In addition, if $\Phi '=e^F \Phi$ is another positive volume form on $M$, then the relation between $\mathcal{M}_\Phi ^\Pi$ and $\mathcal{M}_{\Phi '} ^\Pi$ is
\begin{equation}\label{eq:conformal:volume}
\mathcal{M}_{\Phi '}^\Pi=\mathcal{M}_\Phi ^\Pi - X^\Pi _F.
\end{equation}
Therefore, the modular vector field induces a cohomology class $[\mathcal{M}_\Phi ^\Pi]$ in the first Poisson cohomology group of $M$, which is called  {\it the modular class of $M$}. The Poisson manifold is said to be {\it unimodular} if the modular class vanishes, i.e., $\mathcal{M}_\Phi ^\Pi=X^\Pi _F$, for $F\in C^\infty(M)$. From \eqref{eq:conformal:volume}, it is possible to deduce the following result.
\begin{proposition}\label{rem:preservacion:volumen}
Let $(M,\Pi)$ be a Poisson manifold. Then, $M$ is unimodular if and only if there exists a volume form $\Phi $ which is preserved by all Hamiltonian vector fields. In fact, if $M$ is unimodular and  $\mathcal{M}_{\Phi '}^\Pi=X^\Pi _F$, for a volume form $\Phi '$ and a real function ${F\in C^\infty(M),}$ then the volume form $\Phi = e^{F}\Phi '$ is preserved by all Hamiltonian vector fields.
\end{proposition}

\begin{example}\label{ex:Lie:algebra}
Let $(\mathfrak{g},[\cdot,\cdot]_\g)$ be a Lie algebra of dimension $n$. Then, the modular character $\mathcal{M}_\g\in \mathfrak{g}^*$  of $\mathfrak{g}$ is the trace of the adjoint action
\[
 \mathcal{M}_\g ( \xi )=\mbox{Tr}(\mbox{ad} _\xi ),\mbox{ for $\xi\in {\mathfrak g}.$}
\]
So, if $\{ e_\alpha \}$ is a basis of $\mathfrak{g}$ with 
dual basis $\{e^\alpha \}$, we have
\[
 \mathcal{M}_\g =c^\beta _{\alpha \beta }e^\alpha ,
\]
where $[ e_\alpha ,e_\beta ]_\g =c_{\alpha \beta }^\gamma e_\gamma$.

Now, consider the Lie-Poisson structure $\Pi_{LP}$ on $\mathfrak{g}^*$ given by (\ref{Lie-Poisson}). If $(x_\gamma)$ are global coordinates on ${\mathfrak g}^*$ induced by the basis $\{e^\gamma\}$ of ${\mathfrak g}^*,$  it follows that 
$$\Pi_{LP}=\frac{1}{2}c_{\alpha\beta}^\gamma x_\gamma \frac{\partial }{\partial x_\alpha}\wedge\frac{\partial }{\partial x_\beta}.$$
So, the modular vector field $\mathcal{M}^{\Pi_{LP}} _\Phi$ with respect to the volume form $\Phi=\dif x_1\wedge \ldots \wedge \dif x_n$ is 
\[
\mathcal{M}^{\Pi_{LP}} _\Phi = c^\beta _{\alpha \beta }\frac{\partial}{\partial x_\alpha }.
\]
In fact, following \cite{M}, one may prove that  the linear Poisson structure $\Pi_{LP}$  on ${\mathfrak g}^*$ is unimodular if and only if the Lie algebra ${\mathfrak g}$ is unimodular (that is, the modular character of ${\mathfrak g}$ is zero). 
\end{example}

\begin{example}
Let $(G,\Pi)$ be a connected Poisson-Lie group of dimension $n$ with Lie bialgebra $(\mathfrak{g},\mathfrak{g}^*)$. If $ \mathcal{M}_\g\in \mathfrak{g}^*$  (resp. $ \mathcal{M}_{\g ^*}\in \mathfrak{g}$) is the modular character of $\mathfrak{g}$ (resp. $\mathfrak{g}^*$), then the modular vector field associated with a left-invariant volume form $\nu ^l$ on $G$ is
\begin{equation}\label{eq:mod:vf:Poisson-Lie}
\mathcal{M}^{\Pi} _{\nu ^l} =\frac{1}{2} (\mathcal{M}_{\g^*}^l +  \mathcal{M}_{\g^*}^r +
\Pi^\sharp (\mathcal{M}_{\g}^r))
\end{equation}
(see \cite{ELW}). 
Note that if $G$ is abelian and isomorphic to a vector space $V$, then $\Pi $ is multiplicative if it is linear. Thus, since $V^*$ is a Lie algebra and the Lie algebra in $V$ is trivial, we recover Example \ref{ex:Lie:algebra}.

Now, let $f_0:G\to \R$ be the real $C^\infty$-function on $G$ given by
\begin{equation}\label{f0}
f_0(g) := \mbox{det }(Ad^n_g), \mbox{ for } g\in G,
\end{equation}
where $Ad^n:G\times \wedge^n\g\to \wedge ^n \g$ is the adjoint action of $G$ on $\wedge^n\g.$ 
{Note that $f_0(e) = 1$ and, thus,  since $f_0$ is not zero at every point of $G$ and $G$ is connected, we deduce that $f_0(g) > 0$, for every $g \in G$. Moreover,} {$\mathcal{M}_{\g}^r = \dif (\log f_0)$ (see \cite{ELW})} and 
\[
[\mathcal{M}^{\Pi} _{\nu ^l} ]=  [\frac{1}{2} (\mathcal{M}_{\g^*}^l +  \mathcal{M}_{\g^*}^r+\Pi^\sharp (\dif (\log f_0)))]=[\frac{1}{2} (\mathcal{M}_{\g^*}^l +  \mathcal{M}_{\g^*}^r)].
\]
{So, using that $\Pi(e) = 0$, it is easy to prove that the Poisson manifold $(G, \Pi)$ is unimodular if and only if $\mathcal{M}_{\g^*} = 0$, that is, the dual Lie algebra $\g^*$ is unimodular. Moreover, if $(G, \Pi)$ is unimodular it follows that the modular vector field $\mathcal{M}^{\Pi} _{\nu ^l}$ is just the Hamiltonian vector field of the real function $\frac{1}{2}\log f_0$, that is,}
\[
\mathcal{M}^{\Pi} _{\nu ^l}  = \frac{1}{2} X^{\Pi}_{(\log f_0)}.
\]

\end{example}

\section{Unimodularity of Poisson-Lie groups and preservation of volume forms}

Let $(G,\Pi)$ be a Poisson-Lie group of dimension $n$, $H:G\to \R$ a Hamiltonian function and $\nu\in \wedge^n\g^*,$ with $\nu\not=0.$ Then, a (positive) volume form on $G$ is given by 
$$\Phi=e^\sigma\nu^l,\mbox{ with } \sigma\in C^\infty(G).$$

Moreover, using (\ref{eq:modular:vf}), (\ref{eq:divergencia}) and (\ref{eq:mod:vf:Poisson-Lie}), we deduce that 
\begin{equation}\label{diver-Ham}
\begin{array}{rcl}
{\mathcal L}_{X_H^\Pi}\Phi&=&e^\sigma(X_H^\Pi(\sigma) + {\mathcal M}_{\nu^l}^\Pi (H))\nu^l\\[5pt]&=& e^\sigma(X_H^\Pi(\sigma-log\sqrt{f_0}) + \frac{1}{2}({\mathcal M}_{\g^*}^l(H) + {\mathcal M}^r_{\g^*}(H)))\nu^l,
\end{array}
\end{equation}
with $f_0:G\to \R$ the real $C^\infty$-function on $G$ defined in (\ref{f0}). 

So, from (\ref{diver-Ham}), it follows that 

\begin{theorem}\label{th3.0} Let $(G,\Pi)$ be a connected Poisson-Lie group and $H:G\to \R$ a Hamiltonian function. Then, the Hamiltonian vector field $X_H^\Pi$ preserves a volume form on $G$ if and only if there exists $\sigma\in C^\infty(G)$ such that
$$X_H^\Pi(\sigma-log\sqrt{f_0}) +  \frac{1}{2}({\mathcal M}_{\g^*}^l(H) + {\mathcal M}^r_{\g^*}(H))=0.$$
Moreover, in such a case, if $\nu\in\wedge^n\g^*,$ with $\nu\not=0$, then the volume form $e^\sigma\nu^l$ is preserved by $X_H^\Pi.$ 
\end{theorem}

\begin{remark}\label{r3.1'}
If $g\in G$ is a singular point of $X_{H}^\Pi$ and $X_{H}^\Pi$ preserves a volume form on $G$ then, from Theorem \ref{th3.0}, we deduce that 
$${\mathcal M}^l_{\g^*}(g)(H)+ {\mathcal M}^r_{\g^*}(g)(H)=0.$$

\end{remark}

\begin{remark}\label{r3.2'}
If $H\in C^{\infty}(G)$ is a first integral of vector field ${\mathcal M}_{\mathfrak g^*}^l + {\mathcal M}^r_{\mathfrak g^*}$, then the volume form  $\sqrt{f_0}v^l$ is preserved by the Hamiltonian vector field  $X_H^{\Pi}$.
\end{remark}

Using Theorem \ref{th3.0}, we also obtain the following result 
\begin{corollary}\label{first-theorem}
Let $(G, \Pi)$ be a connected Poisson-Lie group of dimension $n$ such that the dual Lie algebra $\g^*$ is unimodular (that is, the Poisson-Lie structure $\Pi$ is unimodular). If $\nu \in \wedge^n \g^*$ with $\nu \neq 0$, then the volume form $\sqrt{f_0} \nu^l$ is preserved by all Hamiltonian vector fields, where $f_0: G \to \R$ is the real function on $G$ given in (\ref{f0}). 
\end{corollary}

Next, we will prove a converse of this result for a particular kind of Hamiltonian functions: the Morse functions at the identity element of the Poisson-Lie group. 

\begin{definition}
Let $G$ be a Lie group with identity element $e$. A function $H\in C^\infty(G)$ is said to be Morse at $e$ if 
\begin{itemize}
\item[i)] $e$ is a singular point of $H$, that is, $\dif H (e)=0$.
\item[ii)] $(\Hs H )(e)$ is nondegenerate, with $(\Hs H )(e) \colon \g \times \g \to \R$  the
Hessian of $H$ at $e$, i.e., the symmetric bilinear form on $\g$ given by
\begin{equation}\label{eq:Hessian}
(\Hs H )(e) (\xi ,\eta )= \xi (X^{\eta }(H) ), \qquad \xi,\eta \in \g,
\end{equation}
where $X^{\eta}\in \mathfrak{X}(G)$ is an arbitrary vector field such that $X^{\eta}(e)=\eta$.
\end{itemize}
\end{definition}

\begin{remark}
\begin{enumerate}
\item[i)]  We recall that a smooth function $H:G\to \R$ is Morse if the Hessian of $H$ at each singular point of $H$ is nondegenerate. More details on nondegenerate singular points and the definition of the Hessian can be found, for instance, in \cite[Part I]{Mil}).
\item[ii)]  A smooth function $H\colon G\to \R$ is a contrast function on $G$ if $H(e)=0$ and $\dif H(e)=0.$ This kind of functions plays an important role in information geometry (see \cite{Gra}). 
Since the Hamiltonian dynamics induced by $H:G\to \R$ coincides with that induced by the function $H-H(e)$, we can assume (without the loss of generality) that for a Morse function at $e$, $H:G\to \R,$ we have that  $H(e)=0.$ So, a Morse function at $e$ is a contrast function on $G$. \end{enumerate}
\end{remark}

Now, we will prove the result announced above. 
\begin{theorem}\label{th3-2}
Let $(G,\Pi)$ be a connected Poisson-Lie group of dimension $n$ and $H\in C^\infty(G)$ be a function on $G$ which is Morse at $e$.
If $X^{\Pi}_H$ preserves a volume form $\Phi$ on $G$, then the dual Lie algebra $\g^*$ is unimodular and, therefore, the volume form $\sqrt{f_0} \nu^l$, with $\nu\in \wedge^n\g^*$ and $\nu\not=0,$   is preserved by all Hamiltonian vector fields (in particular, for the Hamiltonian vector field of $H$). 
\end{theorem}

\begin{pf}
Let $0\neq \nu \in \wedge ^n \g^*$ and $\nu ^l\in \Omega ^n(G)$ be the corresponding left-invariant volume form on $G$. Suppose that there is a function $\sigma \in C^\infty (G)$ such 
that the volume form $\Phi=e^\sigma \nu ^l$ is preserved by the Hamiltonian vector field $X^{\Pi}_H$, that is,
\[
\Lie _{X^{\Pi}_H} (e^\sigma \nu ^l )=0.
\]
Then, from (\ref{diver-Ham}), it follows that 
\[
X^{\Pi}_H (\sigma ) + \mathcal{M}^{\Pi} _{\nu ^l}(H)=0,
\]
and, using the skew-symmetry of $\Pi$,
\[
\mathcal{M}^{\Pi} _{\nu ^l}(H)=X^{\Pi}_\sigma (H) 
\]
Thus, for any $\xi \in \g$, 
\begin{equation}\label{key-point-4-2}
(\Hs H)(e)( \xi , \mathcal{M}^{\Pi} _{\nu ^l} (e) )=
\xi ( \mathcal{M}^{\Pi} _{\nu ^l}(H) )  = \xi ( X^{\Pi}_\sigma (H) )=(\Hs H)(e)( \xi ,   X^{\Pi}_\sigma (e)  )=0 \mbox{ for all } \xi\in \g,
\end{equation}
because $X^\Pi_F (e)=0$ for any $F\in C^\infty (G)$. Since $(\Hs H)(e)$ is nondegenerate, $\mathcal{M}^{\Pi} _{\nu ^l} (e)=0$. But, from \eqref{eq:mod:vf:Poisson-Lie}, $\mathcal{M}^{\Pi } _{\nu ^l} (e)= \mathcal{M}_{\g^*}$, the modular character of $\g^*$. So, the dual Lie algebra $\g^*$ is unimodular. This, using Corollary \ref{first-theorem}, ends the proof of the result.
\end{pf}
\begin{remark}
Let $(G,\Pi)$ be a connected Poisson-Lie group and $H\in C^\infty(G)$ be a Hamiltonian  function on $G$   such that the identity element of $G$ is a singular point of $H$. Then, following the proof of Theorem \ref{th3-2}, we deduce an interesting result. Namely, if the Hamiltonian vector field of $H$ preserves a volume form, we have that the modular character of the dual Lie algebra $\g^*$ belongs to the kernel of $(Hess \, H)(e).$
\end{remark}

\section{Examples}
\subsection{Kozlov's unimodularity result}
Let $\g$ be a real Lie algebra of dimension $n$ and $\g^*$ be the dual vector space endowed with the Lie-Poisson structure $\Pi_{LP}$ 
(see Example \ref{ex:Lie:algebra}). 

If $I:\g^*\times \g^*\to \R$ is a symmetric $\R$-bilinear form then we can consider the Hamiltonian function $H_I:\frak \g^*\to \R$ given by 
\[
H_I(\mu) = \frac{1}{2}I(\mu,\mu),\mbox{ for } \mu\in \g^*.
\]
So, if $\{e_\gamma\}$ is a basis of $\g$ with dual basis $\{e^\gamma\}$ for $\g^*$ and $\{x_\gamma\}$  the corresponding global coordinates of $\g^*$, we have that 
\[
H_I(x) = \frac{1}{2}I^{\alpha\beta}x_\alpha x_\beta.
\]
It is clear that the identity element of $\g^*$ as an abelian Poisson-Lie group (that is, the zero vector $0$) is a singular point of $H_I.$ Moreover, it is easy to prove that $H_I$ is Morse at $0$ if and only if $I$ is nondegenerate. So, if $I$ is nondegenerate then, using Corollary \ref{first-theorem} and Theorem \ref{th3-2}, we deduce that the Hamiltonian vector field $X_{H_I}^{\Pi_{LP}}$ preserves a volume form on $\g^*$ if only if $\g$ is unimodular. Moreover, in such a case, the function  $f_0:\g^*\to \R$ given by (\ref{f0}) identically vanishes which implies that the Hamiltonian vector field $X_{H_I}^{\Pi_{LP}}$ preserves the volume form $\Phi=dx_1\wedge \dots \wedge dx_n.$ 

This extends a previous result by Kozlov \cite{K88} for the particular case when $H_I$ is a Hamiltonian function of kinetic type, that is, the symmetric bilinear form $I$ is positive definite.

\subsection{Contrast functions on Lie groups}\label{contrast-functions-Lie-groups}

We recall that a smooth function $H\colon G\to \R$ on a Lie group $G$ with with identity element $e$, is a contrast function if $H(e)=0$ $ and $ $\dif H(e)=0.$ It is obvious that if $\Pi$ is a Poisson-Lie structure on $G$ then $X^\Pi_H=X^\Pi_{H-H(e)}$ and so the functions $H$ and $H-H(e)$ define the same dynamics. So, we can always consider the contrast function $H-H(e)$, when $H$ is a Morse function at $e$. This kind of functions plays an important role in information geometry.

Now, assume that $\iota:G\to \R^n$ is an embedding of the Lie group $G$ in $\R^n$, with $\iota(e)=0$. If $\langle\cdot ,\cdot \rangle$ is a non-degenerate bilinear symmetric form on $\R^n$, we can consider the Hamiltonian function $H:G\to \R$ defined by 
$$H=H_{\langle\cdot ,\cdot \rangle}\circ \iota:G\to \R,$$ 
with $H_{\langle\cdot ,\cdot \rangle}:\R^n\to \R$ given by 
$$H_{\langle\cdot ,\cdot \rangle}(x)=\frac{1}{2} \langle x ,x \rangle,\mbox{ for }x\in \R^n.$$

Then, it is clear that the identity element $e$ is a singular point of $H$
and $\dif H(e)=0$, so $H$ is a contrast function. In addition, one can show that if the restriction of the symmetric bilinear form $\langle\cdot ,\cdot \rangle$ to the subspace  $\g=T_eG\subseteq T_0\R^n\cong \R^n$ is nondegenerate  then $e$ is a non-degenerate singular point for $H$. In particular, if $\langle \cdot ,\cdot \rangle$ is a definite positive  symmetric bilinear form on $\R^n$ then $H$ is Morse at $e$.  Thus, if $G$ is a matrix Lie group (a Lie subgroup of $\mathrm{GL}(n,\R)$) then the smooth function $H:G\to \R$ given by
\[
H(A)= \mbox{Tr}\,((Id-A)(Id-A^t)), \qquad A\in G,
\]
is Morse at $Id\in G.$  Here $Id$ is the identity matrix in $\mathrm{GL}(n,\R).$ 

So, from Corollary \ref{first-theorem} and Theorem \ref{th3-2}, it follows that, for a Poisson-Lie structure $\Pi$ on $G,$ the Hamiltonian vector field of $H$ preserves a volume form on $G$ if and only if $\Pi$ is unimodular. 

\begin{remark}
In the particular case when $G=\mathrm{GL}(n,\R)$, the previous function $H$ was used in \cite{Gra}
as a metric contrast function on $\mathrm{GL}(n,\mathbb R)$. We remark again that this kind of function plays an important role in information geometry (see \cite{Gra}, for more details). 
\end{remark}
\subsection{A Hamiltonian system on the Poisson-Lie group $\mathrm{SL}(2,\mathbb R)$}\label{4.3.}
Let $\mathrm{SL}(2,\mathbb R)$ be the special linear group of order $2$, that is, the Lie subgroup of $\mathrm{GL}(2,\R)$ defined by 
\begin{equation}\label{eq:sl2amb}
\mathrm{SL(2,\mathbb R)} = \left\{ A =
\begin{pmatrix}
a_{11} && a_{12} \\
a_{21} && a_{22} \\
\end{pmatrix}
\in \mathrm{GL}(2,\mathbb R) \,|\, \det A = 1
\right\}.
\end{equation}

The vector space $\mathfrak{gl}(2,\R)$ of matrices of order $2$ with real coefficients is the Lie algebra of $\mathrm{GL}(2,\mathbb R)$ and 
\[
\mathfrak{sl}(2,\R)=\{ A\in \mathfrak{gl}(2,\R)\, |\mbox{ Tr}\,A=0\}
\]
is the Lie algebra of $\mathrm{SL}(2,\R).$ 

We can consider the basis $\mathcal B=\{ J_3, J_+,J_-,\displaystyle\frac{1}{2}Id\}$ of $\mathfrak{gl}(2,\mathbb{R})$,  
given by
\[
J_3 = 
\begin{pmatrix}
1 && 0 \\
0 && -1 \\
\end{pmatrix}, \qquad 
J_+ = 
\begin{pmatrix}
0 && 1 \\
0 && 0 \\
\end{pmatrix}, \qquad 
J_- = 
\begin{pmatrix}
0 && 0 \\
1 && 0 \\
\end{pmatrix}. 
\]
Note that $\{J_3,J_{+},J_{-}\}$ is a basis of ${\mathfrak g}=\mathfrak{sl}(2,\R)$ and, moreover, the non-zero commutation relations for the basis ${\mathcal B}$ read 

{
\begin{equation}
\label{eq:sl2alg}
[ J_3,J_{-} ]_{\g} =  -2 J_{-}, \qquad[ J_3,J_{+} ]_{\g} =  2 J_{+}, \quad [ J_+,J_- ]_\g  =J_3.
\end{equation}
}
 The basis of left-invariant vector fields on $\mathrm{GL}(2,\R )$ associated with $\mathcal B$ is just 
\begin{eqnarray*}
 {\frac{1}{2}Id^l} &=& \frac{1}{2}\Big ( a_{11} \frac{\partial}{\partial {a_{11}}}+a_{12} \frac{\partial}{\partial {a_{12}}} + a_{21} \frac{\partial}{\partial {a_{21}}} + a_{22} \frac{\partial}{\partial {a_{22}}} \Big ) ,\\
 J_3^l&=&  a_{11} \frac{\partial}{\partial {a_{11}}}- a_{12} \frac{\partial}{\partial {a_{12}}} + a_{21} \frac{\partial}{\partial {a_{21}}} - a_{22} \frac{\partial}{\partial {a_{22}}} , \\
J_+^l &= & a_{11} \frac{\partial}{\partial {a_{12}}} + a_{21} \frac{\partial}{\partial {a_{22}}} ,\\
J_-^l &=&  a_{12} \frac{\partial}{\partial {a_{11}}} + a_{22} \frac{\partial}{\partial {a_{21}}}.
\end{eqnarray*}

The dual basis of $\mathcal B$ is 
$\{(J^3)^l, (J^+)^l, (J^-)^l  ,\mathrm{d} (\det A) \}$, $\{ J^3, J^+, J^- \}$ being the dual basis of $\{J_3,J_+,J_- \}$.

Next, consider the Poisson-Lie structure on $\mathrm{GL} (2,\mathbb R)$ given by
\begin{equation}
\label{eq:poisssl2amb}
\begin{split}
&\{ a_{11}, a_{12} \} = a_{11} a_{12}, \quad 
\{ a_{11}, a_{21} \} = a_{11} a_{21} , \quad
\{ a_{11}, a_{22} \} = 2 a_{12} a_{21} , \\
&\{ a_{12}, a_{21} \} = 0 , \quad 
\{ a_{12}, a_{22} \} = a_{12} a_{22} , \quad
\{ a_{21}, a_{22} \} = a_{21} a_{22} . 
\end{split}
\end{equation}
This is indeed a Poisson structure on $\mathrm{SL} (2,\mathbb R),$ since $\det A = a_{11} a_{22} - a_{12} a_{21}$ is a Casimir for this structure (see, for instance, 
\cite{Reyman}). 
For this Poisson-Lie group structure, we have that
\begin{equation}\label{eq:sing:points:sl2}
 A=\begin{pmatrix} a_{11}& 0\\ 0 &a_{22} \end{pmatrix}\in GL(2,\mathbb R)
\end{equation}
is a singular point of $\Pi$. 

The previous Poisson-Lie structure is the one defined by the so-called (standard) Drinfel'd-Jimbo $r$-matrix  {
\noindent 
and, interestingly, it is the same as the Sklyanin bracket appeared in the singular value decomposition in \cite{Tomei} for Toda-SVD in \cite{CHU}.}
With respect to the basis  $\{J_3,J_+,J_- \}$,  the (standard) Drinfel'd-Jimbo $r$-matrix is given by 
\begin{equation}
\label{eq:rDJ}
r = J_- \wedge J_+ \in \wedge^2 \mathfrak{sl}(2,\R).
\end{equation}
In order to compute the tangent Lie algebra to this Poisson-Lie structure using that, for coboundary structures, we have 
\[
\delta (\xi) = \mathrm{ad}_\xi (r),
\]
and, from \eqref{eq:sl2alg} and \eqref{eq:rDJ},
\[
\delta (J_3) = 0, \qquad \delta (J_\pm) = J_3 \wedge J_\pm .
\]
Therefore, 
\begin{equation}\label{eq:dualalg}
[ J^3,J^\pm ]_{\g^*} = J^\pm, \qquad [ J^+,J^- ]_{\g^*} = 0 .
\end{equation}
So, the modular character $\mathcal M_{\mathfrak g^*}$ of the dual Lie algebra can be computed from \eqref{eq:dualalg} and reads
\[
\mathcal M_{\mathfrak g^*} =2 J_3\neq 0 .
\]
{
This implies that the Poisson manifold $(\mathrm{SL(2,\mathbb R)},\{ \cdot ,\cdot \})$ is not unimodular. In fact, 
 the modular vector field of the Poisson manifold $(\mathrm{SL(2,\mathbb R)},\{ \cdot ,\cdot \})$ is
\begin{equation}
\label{eq:mod_class_amb}
\frac{1}{2} (\mathcal M_{\mathfrak g^*}^r \mathcal + {\mathcal M}_{\mathfrak g^*}^l ) = 2 a_{11}\frac{\partial}{\partial a_{11}} -2 a_{22}\frac{\partial}{\partial a_{22}},
\end{equation}
which once more proves the Poisson structure is not unimodular.  Note that modular character $\mathcal M_{\mathfrak g}=0$.
}

{
\noindent  Now, we consider }the Hamiltonian function $H\in C^\infty (\mathrm{GL}(2,\mathbb{R}))$ given by 
\begin{equation}
\label{eq:Hsl2amb}
H = \frac{1}{2} \Tr (A^T A) = \frac{1}{2} (a_{11}^2 + a_{12}^2 + a_{21}^2 + a_{22}^2) =  \frac{1}{2} \sum_{i,j=1}^2 a_{ij}^2,
\end{equation}
{
\noindent Note that this Hamiltonian function generates the Toda-SVD flow defined by the Poisson bracket (\ref{eq:poisssl2amb}) (see \cite{Tomei}). }

\noindent Since $\mathrm d H = a_{11} \mathrm d a_{11} + a_{12} \mathrm d a_{12} + a_{21} \mathrm d a_{21} + a_{22} \mathrm d a_{22}$, it follows that
\[
\mathrm{d} H (Id) \not = 0.
\]
Therefore, we cannot use Theorem \ref{th3-2} but nevertheless we can  use Remark \ref{r3.1'}.

\noindent Indeed, the Hamiltonian vector field of $H$ is given by
\begin{eqnarray*}\label{eq:campo:Hamiltoniano:sl2}
X_H^{\Pi} & = & \left( -a_{11} \left(a_{12}^2+a_{21}^2\right)-2(a_{22}-1)a_{12}a_{21} \right) \frac{\partial}{\partial a_{11}} + 
\left(( a_{11}-1)a_{11}a_{22}-( a_{22}-1)a_{12} a_{22} \right) \frac{\partial}{\partial a_{12}} \\ 
& &+ 
\left( (a_{11}-1)a_{11} a_{12}-(a_{22}-1)a_{12} a_{22} \right) \frac{\partial}{\partial a_{21}} + 
\left( 2( a_{11}-1) a_{12} a_{21}+a_{22} \left(a_{12}^2+a_{21}^2\right) \right) \frac{\partial}{\partial a_{22}}\\
&=&X^{11}\frac{\partial}{\partial a_{11}} + X^{12}\frac{\partial}{\partial a_{12}}+X^{21}\frac{\partial}{\partial a_{21}}+X^{22}\frac{\partial}{\partial a_{22}},
\end{eqnarray*}

{\noindent which is the flow of Toda-SVD in dimension two.}
In addition, if $a\not=0$ then, from \eqref{eq:sing:points:sl2}, 
\[
A=\left ( \begin{array}{cc} a & 0 \\0 &\displaystyle\frac{1}{a} \end{array}\right )\in \mathrm{SL}(2,\mathbb R)
\]
is an equilibrium point for the Hamiltonian vector field $X_H^{\Pi}$ on $\mathrm{SL}(2,\R )$. 
 Moreover, using  (\ref {eq:mod_class_amb}) and (\ref{eq:Hsl2amb}), we deduce that 
$${\mathcal M}^l_{\g^*}(A)(H)+{\mathcal M}^r_{\g^*}(A)(H)=\frac{4}{a^2}(a^2-1)(a^2-a+1).$$

In particular, if $a\not=\pm 1$ it follows that 
$${\mathcal M}^l_{\g^*}(A)(H)+{\mathcal M}^r_{\g^*}(A)(H)\not=0,$$ which, from Remark \ref{r3.1'}, implies that $X_H^\Pi$ does not preserve a global volume form on $\mathrm{SL}(2,\R )$.

\subsection{A Hamiltonian system on the Poisson-Lie group  $S^3\cong \mathrm{SU}(2,\C)$}

Let $S^3$ be the unit sphere in $\mathbb{R}^4-\{(0,0,0,0)\}$. In $\mathbb{R}^4-\{(0,0,0,0)\}$ one can consider standard global coordinates 
$(x,y,z,t)$ and the Lie group structure with multiplication
\[
(x,y,z,t)(x',y',z',t')=(xx'-yy'-zz'-tt', xy'+yx'-zt'+tz',zx'-ty'+xz'+yt', zy'+tx'+xt'-yz'), 
\]
which is just the Lie group structure identifying the quaternions $\mathbb{H}$ with $\mathbb{R}^4-\{(0,0,0,0)\}$. Note that the identity element is $(1,0,0,0)$ and that $S^3$ is a closed normal Lie subgroup of $\R^4-\{(0,0,0,0)\}$. If $\{ e_1,e_2,e_3,e_4\}$ is the canonical basis in $\mathbb{R}^4-\{(0,0,0,0)\}$, then the associated left invariant vector fields are
\begin{eqnarray*}
e_1^l&=& x\parcial{x}+y\parcial{y}+z\parcial{z}+t\parcial{t}\\
e_2^l&=& -y \parcial{x}+x\parcial{y}-t\parcial{z}+z\parcial{t}\\
e_3^l&=& -z\parcial{x}+t\parcial{y}+x\parcial{z}-y\parcial{t}\\
e_4^l&=& -t\parcial{x}-z\parcial{y}+y\parcial{z}+x\parcial{t}.
\end{eqnarray*}
Then, the non-zero commuting relations for this basis read 
\[
[e_2,e_3]_\g=-[e_3,e_2]_\g=-2e_4, \qquad [e_2,e_4]_\g=-[e_4,e_2]_\g=2e_3, \qquad [e_3,e_4]_\g=[e_4,e_3]_\g=-2e_2.
\]
In fact, $\{e_2,e_3,e_4\}$ is a basis of the Lie algebra of the Lie subgroup $S^3\cong \mathrm{SU}(2,\C )$. On the other hand, there is a Poisson structure $\Pi$ on $\mathbb{R}^4-\{(0,0,0,0)\}$ whose Poisson bracket $\{\cdot,\cdot\}$ is characterized by the following conditions 
\begin{eqnarray*}
\{ x,y\}=-(z^2+t^2),\\
\{ x,z \} =yz, \quad \{ x,t \} =yt,\\
\{ y,z \} =-xz, \quad \{ y,t \} =-xt.\\
\end{eqnarray*}
Indeed, $\Pi$ is a Poisson-Lie structure on $\R^4-\{(0,0,0,0)\}$ and 
\[
\Pi  = -(z^2+t^2) \parcial{x}\wedge \parcial{y}+yz\parcial{x}\wedge \parcial{z}+yt\parcial{x}\wedge \parcial{t}-xz\parcial{y}\wedge \parcial{z}-xt\parcial{y}\wedge \parcial{t}.
\]
We remark that 
$$ \|\cdot\|:\R^4-\{(0,0,0,0)\}\to \R,\;\;\; (x,y,z,t)\to \|(x,y,z,t)\|^2=x^2+y^2+z^2+t^2$$
is a Casimir function and, thus, $\Pi$ induces a Poisson-Lie structure on $S^3,$ which we also denote by $\Pi.$ It is the standard Poisson-Lie structure on $S^3\cong \mathrm{SU}(2,\C)$ (see \cite{Reyman}). Moreover, using 
 the fact that the adjoint 1-cocycle $\delta$ satisfies
\[
[\xi^l, \Pi ] =  (\delta \xi)^l, \mbox{ for }\xi\in \mathfrak{su}(2),
\]
we deduce that
$$
\delta (e_1)=\delta (e_2)=0,\quad
\delta (e_3)=-e_2\wedge e_3,\quad \delta (e_4)=-e_2\wedge e_4,
$$
and for  the Lie algebra structure on $\g^*$ the non-zero commutation relations read \begin{eqnarray*}
\,[e^2,e^3]_{\g^*}=-[e^3,e^2]_{\g^*}=-e^3, \\
\,[e^2,e^4]_{\g^*}=-[e^4,e^2]_{\g^*}=-e^4.
\end{eqnarray*}
So, the modular character of the dual of the Lie algebra, $\mathcal{M}_{\mathfrak{g}^*}$, is simply $-2e_2$. Therefore,
the modular vector field of the Poisson manifold $(S^3, \Pi )$ is given by
\begin{equation}\label{eq:su(2):modular:class}
\frac{1}{2} (  \mathcal{M}_{\mathfrak{g}^*}^l+   \mathcal{M}_{\mathfrak{g}^*}^r)= -2(x\parcial{y}-y\parcial{x} ).
\end{equation}

Thus,  it follows that $(S^3,\{ \cdot ,\cdot \})$ is not unimodular.

Now, it is easy to show that for any smooth functions $F\colon \mathbb{R}\to \mathbb{R}$ and $G\colon \mathbb{R}^2\to \mathbb{R}$, the function
\[
\begin{array}{rccl}
H\colon &\mathbb{R}^4-\{(0,0,0,0)\}& \to & \mathbb{R} \\
              & (x,y,z,t) & \mapsto & F(x^2+y^2)+G(z,t)
\end{array}
\]
is a smooth first integral of  $\frac{1}{2} (  \mathcal{M}_{\mathfrak{g}^*}^l+   \mathcal{M}_{\mathfrak{g}^*}^r)$. In particular, the Hamiltonian function $$H\colon S^3\to \R,\quad (x,y,z,t)\mapsto P(z,t),$$ where $P$ is an arbitrary smooth function on $\mathbb{R}^2$, is a first integral of 
$\frac{1}{2} (  \mathcal{M}_{\mathfrak{g}^*}^l+   \mathcal{M}_{\mathfrak{g}^*}^r)_{|S^3}$. Therefore, using Theorem \ref{th3.0}, we deduce that the Hamiltonian vector field $X^{\Pi}_H$ preserves any left-invariant volume form on $S^3$.

\bigskip

{

%
%

}
\subsection{A deformation of a conservative Lorenz system}

In \cite{BMR2017biHamiltonian}, the authors  presented a method for obtaining integrable deformations of Lie-Poisson bi-Hamiltonian systems. The initial data is a dynamical system $\mathcal D$ on $\mathbb{R}^n$, which is bi-Hamiltonian with respect to two compatible linear Poisson structures. Then, they construct a bi-Hamiltonian deformation $\mathcal D_\eta$ of the dynamical system $\mathcal D$, whose bi-Hamiltonian structure is provided by a pair of Poisson-Lie structures on a non-abelian Lie group $G_\eta$. Under certain conditions, the limit $\eta\to 0$ of $\mathcal D_\eta$ is just the initial dynamical system $\mathcal D$. The Poisson coalgebra approach to Hamiltonian systems was discussed in \cite{BaRa}.

In what follows,  we illustrate our results with two examples coming from the Poisson-Lie deformation theory of Lie-Poisson bi-Hamiltonian  systems. 

As a first example, consider the integrable deformation of a conservative Lorenz system studied in \cite{BBM2016rosslerlorentz}, whose equations of motion read
\begin{equation}
\begin{array}{rcl}
\dot x &=& y, \\
\dot y &=&  \displaystyle\frac{x}{2} \left( 4 + \displaystyle\frac{e^{-2 \p (z + w)} - 1}{\p} + \p (2 x^2 - y^2) \right), \\
\dot z &=& x y , \\
\dot w &=& 0.
\end{array}
\label{eq:deformedLorenz}
\end{equation}
This dynamical system can be written as a Hamiltonian system on the 4-dimensional Lie group $G_{\eta} \simeq B_2 \times B_2$, where $B_2$ is the unique connected simply connected non-abelian Lie group of dimension $2$. We start by introducing exponential coordinates of the second kind on $G_{\eta}$, defined by
\begin{equation}\label{eq:exp:coord}
g (x,y,z,w) = \exp( x \, X) \exp( y \, Y) \exp( z \, Z) \exp( w \, W) .
\end{equation}
In terms of these coordinates the group law reads
{\small
\begin{equation}
\begin{array}{rcl}
g (x,y,z,w) \cdot g '(x',y',z',w') &=&
 g \bigg(x+ \displaystyle\frac{e^{ -\p (z + w)}}{2} \bigg( 2 x' \cosh \left( \displaystyle\frac{\p w}{\sqrt 2} \right) - \sqrt 2 y' \sinh \left(  \displaystyle\frac{\p w}{\sqrt 2} \right)  \bigg),
\\[8pt]&&y+e^{ -\p (z + w)} \bigg( y' \cosh \left(  \displaystyle\frac{\p w}{\sqrt 2} \right) - \sqrt 2 x' \sinh \left(  \displaystyle\frac{\p w}{\sqrt 2} \right)  \bigg)
, z+z', w+w' \bigg) 
\label{eq:glawLorenz}
\end{array}
\end{equation}
}
where $\eta > 0$ and the Poisson bivector field is given by
\begin{equation}
\Pi_{B_2 \times B_2} = \frac{e^{-2 \p (z + w)} - 1 + \p^2 (2 x^2 - y^2)}{4 \p} \frac{\partial}{\partial x} \wedge \frac{\partial}{\partial y} + 
\frac{y}{2} \frac{\partial}{\partial x} \wedge \frac{\partial}{\partial z} + x \frac{\partial}{\partial y} \wedge \frac{\partial}{\partial z} .
\label{eq:piLorentz}
\end{equation}
$\Pi_{B_2 \times B_2}$ defines a Poisson-Lie structure on $G_\eta$ and the dynamical system \eqref{eq:deformedLorenz} can be recovered by means of the Hamiltonian function 
\begin{equation}
\mathcal H = 2 (z-w) - x^2 .
\label{eq:HLorenz}
\end{equation}
The Lie algebra $\mathfrak g_{\eta}$ associated to the Lie group $G_{\eta}$ is characterized by the commutators 
\begin{equation}
[X,Z]_\eta=\eta X, \quad [X,W]_\eta=\eta (X +Y),\quad [Y,Z]_\eta=\eta Y, \quad [Y,W]_\eta=\eta \left( \frac{1}{2} X + Y \right), \quad [X,Y]_\eta=[Z,W]_\eta=0.
\end{equation}
To prove that the Poisson-Lie structure \eqref{eq:piLorentz} is unimodular and therefore the system \eqref{eq:deformedLorenz} admits invariant volume forms, we need to compute the commutator map $\delta_{B_2 \times B_2} = d_e \Pi_{B_2 \times B_2} : {\mathfrak g}\to \wedge^2 {\mathfrak g}$. Explicitly, we have
\begin{equation}
\delta_{B_2 \times B_2} (X) = Y \wedge Z, \qquad \delta_{B_2 \times B_2} (Y) = \frac{1}{2} X \wedge Z, \qquad \delta_{B_2 \times B_2} (Z) = - \frac{1}{2} X \wedge Y, \qquad \delta_{B_2 \times B_2} (W) = - \frac{1}{2} X \wedge Y .
\end{equation}
From this cocommutator, we directly obtain that the Lie algebra $\mathfrak g^*$ of the dual Poisson-Lie group is 
\begin{equation}
[\bar X, \bar Y]_{\mathfrak g^*} = -\frac{1}{2} (\bar Z + \bar W), \qquad [\bar X, \bar Z]_{\mathfrak g^*} = \frac{1}{2} \bar Y, \qquad [\bar Y, \bar Z]_{\mathfrak g^*} = \bar X, \qquad [\bar W, \bar X]_{\mathfrak g^*} = [\bar W, \bar Y]_{\mathfrak g^*} = [\bar W, \bar Z]_{\mathfrak g^*} = 0 .
\end{equation}
Therefore, $\mathfrak g^* \simeq \mathfrak{sl}(2, \mathbb R) \oplus \mathbb R$, which is clearly unimodular since it is a direct sum of two unimodular Lie algebras. 

\noindent To have an invariant volume form, we obtain the left invariant vector fields as
\begin{equation}
\begin{array}{rcl}
X^l&= & e^{ -\p (z + w)} \cosh( \displaystyle\frac{\eta w}{\sqrt 2})\frac{\partial}{\partial x}-\sqrt 2e^{ -\p (z + w)}  \sinh( \displaystyle\frac{\eta w}{\sqrt 2})\frac{\partial}{\partial y}\\[6pt]

Y^l&=& - \displaystyle\frac{\sqrt 2}{2}e^{ -\p (z + w)}  \sinh( \displaystyle\frac{\eta w}{\sqrt 2}) \displaystyle\frac{\partial}{\partial x}+ e^{ -\p (z + w)}  \cosh \displaystyle(\frac{\eta w}{\sqrt 2}) \displaystyle\frac{\partial}{\partial y}\\[6pt]
Z^l&=& \displaystyle\frac{\partial}{\partial z}\\[6pt]
W^l&=& \displaystyle\frac{\partial}{\partial w}\\
\end{array}
\end{equation}
and the right invariant vector fields as
\begin{equation}
\begin{array}{rcl}
X^r&= & \displaystyle \frac{\partial}{\partial x}\\[6pt]

Y^r&=&  \displaystyle\frac{\partial}{\partial y}\\[6pt]
Z^r&=&-\p  x \displaystyle\frac{\partial}{\partial x}-\p  y\frac{\partial}{\partial y}+ \displaystyle\frac{\partial}{\partial z}  \\[6pt]
W^r&=&-(\p  x+  \displaystyle\frac{\p}{2} y)\frac{\partial}{\partial x}- (\p  x+\p y) \displaystyle\frac{\partial}{\partial y}+ \displaystyle\frac{\partial}{\partial w}.\\
\end{array}
\end{equation}
\noindent For invariant $1$-forms we have that 
\begin{equation}
\begin{array}{rcl}
 \bar X^l&= &  e^{ \p (z + w)} \cosh( \displaystyle\frac{\eta w}{\sqrt 2})\dif  x+\frac{\sqrt 2}{2} e^{ \p (z + w)}  \sinh( \displaystyle\frac{\eta w}{\sqrt 2})\dif y\\[10pt]
\bar Y^l&=& \sqrt 2 e^{ \p (z + w)}  \sinh(\frac{\eta w}{\sqrt 2})\dif  x+ e^{ \p (z + w)}  \cosh( \displaystyle\frac{\eta w}{\sqrt 2})\dif  y\\[6pt]
\bar Z^l&=& \dif  z\\[6pt]
\bar W^l&=&\dif  w\\
\end{array}
\end{equation}
and 
\begin{equation}
\begin{array}{rcl}
\bar X^r&= & \dif  x-\p x \dif z +(\p x+\p y)\dif w\\[6pt]

\bar Y^r&=& \dif  y-\p y \dif  z+(\p x +\p y)\dif  w\\[6pt]
\bar Z^r&=&\dif  z \\[6pt]
\bar W^r&=&\dif  w,\\
\end{array}
\end{equation}
\noindent The modular character of Lie algebra $\g_{\p}$ is $\mathcal{M}_{\g_{\eta}}=-2\eta(\bar Z+\bar W)$ and so $$\mathcal{M}_{\g_{\eta}}^r=-2\eta(\bar Z^r+\bar W^r)=-2\eta(\dif z+\dif w)=\dif (\log f_0).$$ Therefore, $f_0=e^{-2\p (z+w)}$. For $v=\bar X\wedge \bar Y\wedge \bar Z\wedge\bar  W\in \Lambda^4\mathfrak g^*$, we get the invariant volume form
\begin{equation}
\Phi=\sqrt f_0 v^l=e^{ \p (z+w)}\dif x\wedge \dif y\wedge \dif z\wedge \dif w.
\end{equation}
 { It is easy to see that $\mathcal L_{{X_{\mathcal H}^\Pi}} \Phi = \dif (i_{{X_{\mathcal H}^\Pi} }\Phi) =0$. }
 

In fact, it is easy to show that $\mathcal L_{X_{\mathcal F}^{\Pi}} \Phi = 0$, for an arbitrary smooth function $\mathcal F : G_{\eta} \to \mathbb R$, in agreement with Corollary \ref{first-theorem}. 


\begin{remark}
As stated above, \eqref{eq:deformedLorenz} is an integrable deformation of a different system, which we call a conservative Lorenz system (see \cite{BBM2016rosslerlorentz} for details), with equations
\begin{equation}
\dot x = y, \qquad \dot y = x(1-z), \qquad \dot z = x y .
\label{eq:Lorenz}
\end{equation}
Since $B_2 \times B_2$ is diffeomorphic to $\mathbb R^4$ as a smooth manifold, this implies that we can recover the previous dynamical system \eqref{eq:Lorenz} starting from \eqref{eq:deformedLorenz} and taking the $\p \to 0$ limit. Moreover, we can also take the $\p \to 0$ limit in the group law \eqref{eq:glawLorenz} obtaining the abelian Lie group structure on $\mathbb R^4$
\begin{equation}
g (x,y,z,w) \cdot g '(x',y',z',w') = g \big(x + x', y + y', z+z', w+w' \big) 
\label{eq:glawLorenzabelian}
\end{equation}
and in the Poisson-Lie bivector 
\begin{equation}
\Pi_{\mathbb R^4} = -\frac{1}{2} (x+w) \frac{\partial}{\partial x} \wedge \frac{\partial}{\partial y} + 
\frac{y}{2} \frac{\partial}{\partial x} \wedge \frac{\partial}{\partial z} + x \frac{\partial}{\partial y} \wedge \frac{\partial}{\partial z} ,
\label{eq:piLorenzund}
\end{equation}
if we take the $\eta\to 0$ limit in \eqref{eq:piLorentz}.

Note that this Poisson-Lie structure is nothing more than the Lie-Poisson structure associated with $\mathfrak g^*$, i.e. the linear Poisson structure on the vector space $(\mathfrak g^*)^* \simeq \mathfrak g$ associated with the Lie algebra $\mathfrak g^* \simeq \mathfrak{sl}(2, \mathbb R) \oplus \mathbb R$. Therefore, using the Hamiltonian \eqref{eq:HLorenz} we obtain the system 
\begin{equation}
\dot x = y, \qquad \dot y = x(2-z-w), \qquad \dot z = x y, \qquad \dot w = 0 .
\end{equation}
From this system, in order to recover \eqref{eq:Lorenz}, we simply restrict to the submanifold $\mathcal{S} = \{ (x,y,z,w) \in \mathbb R^4 | w=1\}$, taking into account that \eqref{eq:piLorenzund} is tangent to $\mathcal S$.

\end{remark}

\subsection{A deformed Euler top}

Following a similar procedure as above we can consider the integrable deformation of the Euler top (see \cite{BMR2017biHamiltonian} for more details), given by
\begin{equation}
\begin{split}
\dot x &= e^{\p x} (y^2 - z^2), \\
\dot y &= \p e^{\p x} y z^2 - \frac{1}{2} \p e^{\p x}  y (y^2+z^2) + \frac{\sinh(\p x) }{\p} (2 z - y), \\
\dot z &= -\p e^{\p x} y^2 z + \frac{1}{2} \p e^{\p x} z (y^2 + z^2)  + \frac{\sinh(\p x)}{\p} (z - 2 y) .
\end{split}
\label{eq:deformedtop}
\end{equation}
This system can be written as a bi-Hamiltonian system on the so-called ``book'' Lie group $G_{\eta}$, with the Lie algebra $\mathfrak g$ defined by
\[
[X,Y]=-\eta Y, \quad [X,Z]=-\eta Z, \quad [Y,Z]= 0,
\]
 where $\eta\in\mathbb R$ is the deformation parameter. Since $G_{\eta}$ is diffeomorphic to $\mathbb R^3$ we can choose global coordinates such that the group law reads
\begin{equation}
g (x,y,z) \cdot g '(x',y',z') = g \big(x + x', y + y' e^{- \p x}, z + z' e^{- \p x} \big) .
\label{eq:glawbook}
\end{equation}

The system \eqref{eq:deformedtop} is bi-Hamiltonian since we can find two different Poisson-Lie structures on $G_{\eta}$, namely
\begin{equation}
\Pi_{G,0} = - z \frac{\partial}{\partial x} \wedge \frac{\partial}{\partial y} + y \frac{\partial}{\partial x} \wedge \frac{\partial}{\partial z} + \frac{1}{2} \bigg( - \p (y^2 + z^2) + \frac{e^{- 2 \p x} - 1}{\p} \bigg) \frac{\partial}{\partial y} \wedge \frac{\partial}{\partial z}
\label{eq:pitop0}
\end{equation}
and 
\begin{equation}
\Pi_{G,1} = - y \frac{\partial}{\partial x} \wedge \frac{\partial}{\partial y} + z \frac{\partial}{\partial x} \wedge \frac{\partial}{\partial z} + \bigg( - \p y z + \frac{e^{- 2 \p x} - 1}{\p} \bigg) \frac{\partial}{\partial y} \wedge \frac{\partial}{\partial z} ,
\label{eq:pitop1}
\end{equation}
and two different functions 
\begin{equation}
\mathcal H_0 = y z e^{\p x} + 2 \left( \frac{\cosh(\p x) - 1}{\p^2} \right),
\label{eq:Htop0}
\end{equation}
and
\begin{equation}
\mathcal H_1 = -\frac{1}{2}  (y^2 + z^2) e^{\p x} +  \frac{\cosh(\p x) - 1}{\p^2}  ,
\label{eq:Htop1}
\end{equation}
such that \eqref{eq:deformedtop} is the dynamical system defined by the Hamiltonian vector field $\Pi_{G,0}^\sharp (d \mathcal H_0) = \Pi_{G,1}^\sharp (d \mathcal H_1)$. Note that both Poisson-Lie structures form a Poisson pencil on $G_{\eta}$, \textit{i.e.} the convex combination $\lambda \Pi_0 + (1 - \lambda) \Pi_1$ is also a Poisson structure, for $\lambda\in \R.$

In this case, the commutator maps $\delta_{G,i} = d_e \Pi_{G,i}$ associated to \eqref{eq:pitop0} and \eqref{eq:pitop1}, respectively, read 
\begin{equation}
\delta_{G,0} (X) = - Y \wedge Z, \qquad \delta_{G,0} (Y) = X \wedge Z, \qquad \delta_{G,0} (Z) = - X \wedge Y ,
\end{equation}
and
\begin{equation}
\delta_{G,1} (X) = - 2 Y \wedge Z, \qquad \delta_{G,1} (Y) = - X \wedge Y, \qquad \delta_{G,1} (Z) = X \wedge Z .
\end{equation}
and therefore the dual Lie algebras read
{ \begin{equation}
\begin{split}
&[\bar X, \bar Y]_{\mathfrak g^*_0} = - \bar Z, \qquad [\bar X, \bar Z]_{\mathfrak g^*_0} = \bar Y, \qquad [\bar Y, \bar Z]_{\mathfrak g^*_0} = - \bar X, \\
&[\bar X, \bar Y]_{\mathfrak g^*_1} = - \bar Y, \qquad [\bar X, \bar Z]_{\mathfrak g^*_1} = \bar Z, \qquad [\bar Y, \bar Z]_{\mathfrak g^*_1} = -2 \bar X, \\
\end{split}
\end{equation}}
{Since $\mathfrak g^*_0 \simeq \mathfrak{so}(3)$ and $\mathfrak g^*_1 \simeq \mathfrak{sl}(2, \mathbb R)$ are both unimodular Lie algebras, both Poisson-Lie structures \eqref{eq:pitop0} and \eqref{eq:pitop1}  are unimodular, and therefore the dynamical  system \eqref{eq:deformedtop} admits invariant volume forms.}

Taking the $\p \to 0$ limit in the equations \eqref{eq:glawbook}, \eqref{eq:pitop0} and \eqref{eq:pitop1} we obtain the 3-dimensional abelian Lie group $\mathbb R^3$ together with the Poisson bivectors
\begin{equation}
\begin{split}
\Pi_{\mathbb R^3,0} &= - z \frac{\partial}{\partial x} \wedge \frac{\partial}{\partial y} + y \frac{\partial}{\partial x} \wedge \frac{\partial}{\partial z} - x  \frac{\partial}{\partial y} \wedge \frac{\partial}{\partial z} , \\
\Pi_{\mathbb R^3,1} &= - y \frac{\partial}{\partial x} \wedge \frac{\partial}{\partial y} + z \frac{\partial}{\partial x} \wedge \frac{\partial}{\partial z} - 2 x  \frac{\partial}{\partial y} \wedge \frac{\partial}{\partial z} ,
\label{eq:pitop0und}
\end{split}
\end{equation}
which are precisely the Lie-Poisson structures associated with the Lie algebras $\mathfrak g^* \simeq \mathfrak{so}(3)$ and $\mathfrak g^* \simeq \mathfrak{sl}(2, \mathbb R)$, respectively. Also, the Hamiltonian functions \eqref{eq:Htop0} and \eqref{eq:Htop1} go to 
\begin{equation}
\mathcal H_0^{\p = 0} = x^2 + y z ,
\label{eq:Htop01}
\end{equation}
and
\begin{equation}
\mathcal H_1^{\p = 0} = -\frac{1}{2} (x^2 + y^2 + z^2) .
\label{eq:Htop11}
\end{equation}
The dynamical system defined by the Hamiltonian vector field $  \Pi_{\mathbb R^3,0}^\sharp (d \mathcal H_0^{\p = 0}) = \Pi_{\mathbb R^3,1}^\sharp (d \mathcal H_1^{\p = 0})$ reads
\begin{equation}
\dot x = y^2 - z^2, \qquad \dot y = x(2z-y), \qquad \dot z = x (z-2y) ,
\label{eq:top}
\end{equation}
which is equivalent to a particular case of the Euler top.

 \noindent Now, we obtain the volume form preserved by the flow of the Hamiltonian vector field. 	A basis of left and right-invariant vector fields for the basis $\{\partial_x,\partial_y,\partial_z\}$ on ${\mathfrak g}_\eta,$ are
 			\[
 			 X^l=\partial_x,\quad  Y^l=e^{-\eta x}\partial_y,\quad   Z^l=e^{-\eta x}\partial_z,
 			\]
 			\[
 			 X^r= \partial_x-\eta y  \partial_y-\eta z  \partial_z,\quad  Y^r= \partial_y,\quad  Z^r= \partial_z,
 			\]
	and a basis of the	left and right invariant one forms for the dual basis $\{dx,dy,dz\}$ on ${\mathfrak g}^*_\eta$, are
 			\[
 			\bar X^l=dx,\quad 	\bar Y^l=e^{\eta x}dy,\quad 	\bar Z^l=e^{\eta x}dz.
 			\]	
 			\[
 				\bar X^r=dx,\quad 	\bar Y^r=dy+\eta y dx,\quad 	\bar Z^r=dz+\eta z dx.
 			\]

\noindent  The modular character of Lie algebra $\mathfrak g$ is $\mathcal{M}_{\g}=-2\eta\bar X$ and so $\mathcal{M}_{\g}^r =-2\eta \dif x= \dif (\log f_0)$. Therefore, $f_0=e^{-2\eta x}$. Then, by the Corollary \ref{first-theorem}, for $v=\bar X\wedge \bar Y\wedge \bar Z\in \wedge^3\mathfrak g^*$, the volume form 
 \[
 \Phi=\sqrt{f_0} \nu^l=e^{\eta x}dx\wedge dy\wedge dz
 \]
 is preserved by all Hamiltonian vector fields.
 

 In fact, it is easy to show that $\mathcal L_{X_{\mathcal G}^{\Pi_{G,0}}} \Phi = L_{X_{\mathcal F}^{\Pi_{G,1}}} \Phi = 0$, for arbitrary functions $\mathcal G, \mathcal F : G_{\eta} \to \mathbb R$, in agreement with Corollary \ref{first-theorem}.


  {\cbl 	

}

\section{Conclusions and future work}
We have discussed the existence of invariant volume forms for Hamiltonian systems on Poisson-Lie groups. In particular, if the identity element of the Lie group is a nondegenerate singularity of the Hamiltonian function, the Hamiltonian system admits an invariant volume form if and only if the Poisson-Lie structure is unimodular. Several examples which illustrate our theoretical results are also exhibited. 

We also present some examples of Hamiltonian functions on Poisson-Lie groups with non-degenerate singularity at the identity element which appear, in a natural way, as metric contrast functions (a special class of divergence functions in information geometry) on matrix groups (see Section \ref{contrast-functions-Lie-groups}).

On the other hand, in \cite{LeZh} the authors discuss a connection between information geometry and discrete geometric mechanics. More precisely, they present the relationship between divergence functions (in information geometry) and (regular) discrete Lagrangian functions (in discrete geometric mechanics). 

So, following the previous ideas, it would be interesting to study the relation between Hamiltonian dynamics on Poisson-Lie groups, discrete geometric mechanics and information geometry. This will be postponed for a future research.  

In another direction, an interesting class of Poisson manifolds closely related to Poisson-Lie groups are Poisson homogeneous spaces, given by the quotient $G/H$ of a Poisson-Lie group $G$ with a closed Lie subgroup $H$ of $G$ (see \cite{STS}). As an extension  of the results in this paper, it would be interesting to discuss the existence of invariant volume forms for Hamiltonian systems on $G/H$. In this setting we remark the following facts: 

\begin{itemize}
\item As in the case of Poisson-Lie groups, the explicit description of the symplectic leaves and the symplectic structure on them for the Poisson structure on $G/H$ is, in general, a hard problem. 

\item The Poisson homogeneous space $G/H$ is integrable and an explicit symplectic groupoid ${\mathcal G}(G/H)$ integrating $G/H$ may be constructed  (see \cite{Lu} and also \cite{BoCiStaTa1, BoCiStaTa2}). But, in general, the dimension of ${\mathcal G}(G/H)$ is $2\dim (G/H)$. 
\end{itemize}

So, in order to study the existence of invariant volume forms, for a Hamiltonian system on $G/H$, it  seems reasonable to discuss the problem directly in the homogeneous Poisson space $G/H$. In fact, a description of the modular class of the Poisson structure on $G/H$ has been given in \cite{Lu} (see also \cite{Caseiro, ELW}) and, therefore, we could discuss the relation between the unimodularity of the Poisson structure and the existence of invariant volume forms for the Hamiltonian system.

\noindent{\bf Acknowledgments:} The authors thanks the anonymous referees for their valuable comments. D. Iglesias, JC Marrero and E. Padr\'on acknowledge financial support from the Spanish Ministry of
Science and Innovation under grant PGC2018-098265-B-C32. I. Guti\'errez-Sagredo acknowledge financial support from Agencia Estatal de Investigaci\'on (Spain) under grant PID2019-106802GB-I00/AEI/10.13039/50110 0011033.

\end{document}